\documentclass[11pt]{amsart}
\usepackage[english]{babel}
\usepackage{amssymb}
\usepackage{amsmath}
\newcommand{\can}{\overline{\phantom{x}}}

\newtheorem{lemma}{Lemma}
\newtheorem{theorem}{Theorem}
\newtheorem{proposition}{Proposition}
\newtheorem{corollary}{Corollary}

\theoremstyle{definition}

\newtheorem{example}{Example}

\newtheorem{remark}{Remark}

\keywords{Forms of higher degree, homogeneous polynomials, multiplicative forms, strongly multiplicative forms.}

\subjclass[2000]{Primary: 11E76; Secondary: 11E04, 12E05}

\title{Some classes of multiplicative forms of higher degree}

\author{S. Pumpl\"un}
\email{susanne.pumpluen@nottingham.ac.uk}
\address{School of Mathematics\\
University of Nottingham\\
University Park\\
Nottingham NG7 2RD\\
United Kingdom
}

\begin{document}

\maketitle

\begin{abstract}  Several notions of multiplicativity are introduced for forms of  degree $d\geq 3$ over a field of characteristic
0 or greater than $d$.
Examples of multiplicative and strongly multiplicative forms of higher degree are given. Conditions
 restricting the structure of a strongly multiplicative form are found.
\end{abstract}

\section*{Introduction}

 Pfister forms play an important role in the algebraic theory of quadratic forms  over fields of characteristic not 2,
 and have a number of remarkable properties. It would be interesting to see whether there are forms of higher degree
 which display a similar behaviour.

 Let $d\geq 2$ be an integer and let $k$ be a field such that ${\rm char}\,k=0$ or $\mbox{char}\,k>d$.
 Let $\varphi:V\to k$ be a  form of degree $d$ on a $k$-vector space
$V$ of dimension $n$ (i.e., after suitable identification,  $\varphi$ is a homogeneous polynomial of degree $d$ in $n$ indeterminates).
Let
 $L= k(X,Y)$ and $K=k(X)$ where $X=(x_1,\dots,x_n)$
and $Y=(y_1,\dots,y_n)$ are sets of independent indeterminates over $k$. The form $\varphi$ is called
\emph{multiplicative} if
$$\varphi(X)\cdot\varphi(Y)=\varphi(z_1,\dots,z_n)$$
for suitable $z_1,\dots,z_n\in L$.

This reminds of forms of higher degree
which permit composition, only that in that case
 the $z_l$ are required to be bilinear forms in the $x_i$'s and $y_i$'s. Thus, as in the $d=2$ case, we may think of
the formula above as some kind of generalized composition.
$\varphi$ is called {\it strongly multiplicative} if
$\varphi(X)\varphi\cong \varphi$ over $K$.

 For $d=2$ the situation is known: a quadratic form $\varphi$  is multiplicative if and only if $D_l(\varphi)$ is a group for every field extension
$l$ over $k$ [L, p.~324]. Every nondegenerate isotropic quadratic form is multiplicative and every
hyperbolic quadratic form is strongly multiplicative.
 The anisotropic
multiplicative or strongly multiplicative quadratic forms can be completely classified: for any anisotropic
nondegenerate quadratic form
$q$ over $k$, $q$ is a Pfister form  if and only if $q$ is multiplicative if and only if $q$ is strongly multiplicative.
Thus a multiplicative quadratic form is either an isotropic form or an anisotropic Pfister form.
A strongly multiplicative form is either a hyperbolic form or an anisotropic Pfister form [L, p.~325].

For $d\geq 3$ it makes sense to also introduce variations of these notions.
  $\varphi$ is called \emph{Jordan multiplicative} if
$$\varphi(X)^2\cdot\varphi(Y)=\varphi(z_1,\dots,z_n)$$
for suitable $z_1,\dots,z_n\in L$, and {\it strongly Jordan multiplicative} if
$$\varphi(X)^2 \varphi\cong \varphi$$
over $K$.

We initiate the study of multiplicative forms of higher degree  in one of the senses just introduced,  focusing mostly on strongly (Jordan)
 multiplicative forms.

Canonical candidates for ``two-fold cubic Pfister forms'' are the reduced norms of central simple associative
algebras of index 3 and for ``three-fold cubic Pfister forms'' the generic norms
of 27-dimensional exceptional simple Jordan algebras (also called {\it Albert algebras}). All these are multiplicative.
For cubic forms, the Tits process (cf. for instance Petersson-Racine [P-S1, 2]) plays a  role similar to that
of the  Cayley-Dickson doubling:
if $J$ is a simple cubic Jordan algebra over a field $k$ of characteristic not 3, then $J$ can be obtained by successive applications of the
 Tits process starting from $k1$.
However, the norm of an algebra obtained via the first Tits construction is determined not only by the norm of
the associative algebra $A$  and the scalar in $k^\times$ one starts out with, but also by the multiplication
in $A$ (more precisely, by the trace of $A$). This is in sharp contrast to the situation for quadratic forms
where only the norm of the algebra and a scalar is needed in order to perform the Cayley-Dickson doubling.
In this light, it will probably not be enough to study for which strongly multiplicative cubic forms $\varphi$ and
scalars $a\in k^\times$ the form $\varphi\perp a\varphi$ is strongly multiplicative. One will have to look at cubic forms
of the kind $\varphi(x)+ a\varphi(y)+ a^2\varphi(z)+\theta(x,y,z)$ with $\theta$ some trilinear form.
It remains unclear for now how to proceed further. Moreover, there might be
a way to generalize the first Tits construction to get 81-dimensional and then 243-dimensional
 cubic forms (etc.) which are strongly multiplicative.
It seems reasonable to expect, however, that Pfister's theory will have no analogue for forms of higher degree.
Nonetheless there is  hope to find interesting phenomena at least for certain classes of  cubic or quartic forms.

Definitions and basic results for multiplicative forms are collected in the first section of the present paper.

Sections 2 and 3  are of a somewhat expository nature since they deal with three
 classes of  forms which are known but have not been treated in our context so far:
 We consider nondegenerate higher degree  forms permitting composition or Jordan composition
which were classified by Schafer [S], and forms belonging to
simple structurable algebras with involution with the property that the space of skew-hermitian
elements has dimension 1. The best-known example of such an algebra is probably
the 56-dimensional irreducible module $\mathcal{M}$ for the split simple Lie algebra
of type $E_7$. Structurable algebras were used to construct all classical simple isotropic Lie algebras [A2].
 Examples of such algebras can be constructed out of separable Jordan algebras
over $k$ of degree 3, cf. Allison [A1], Springer [Sp1, 2]. There exists a Cayley-Dickson process for constructing a class of such algebras by
endowing the direct sum of two copies of a Jordan algebra of degree 4 over $k$ which carries a nondegenerate quartic form
permitting Jordan composition, with the structure of an algebra with involution, see Allison-Faulkner
[A-F1]. It suggests that also for quartic forms
any construction to obtain new strongly Jordan
multiplicative forms out of old ones will probably involve more information than just a scalar and
 the form one starts out with, contrary to the quadratic case.
In section 4, we build new examples of (strongly) multiplicative forms, for instance by using forms of lower degree as
"building blocks". In Section 5, we restrict the structure of strongly multiplicative forms.
It turns out that there are no strongly multiplicative forms of degree $d\geq 3$ which are separable
 other than $\varphi\cong\langle 1\rangle$ and no strongly Jordan multiplicative  separable cubic
 forms other than $\varphi\cong\langle 1\rangle$.
 We close with some observations on other possible concepts of multiplicity for forms of higher degree in Section 6.

\section{Preliminaries}

\subsection{}  Let $k$ be a field of characteristic 0 or $> d$.
A $d$-{\it linear form} over $k$ is a $k$-multilinear map $\theta : V \times
\dots \times V \to k$ ($d$-copies) on a finite-dimensional vector space $V$ over $k$ which is {\it symmetric};
 i.e., $\theta (v_1,
\dots, v_d)$ is invariant under all permutations of its variables.
A {\it form of degree $d$} over $k$ is a map $\varphi:V\to k$  on a finite-dimensional vector space $V$ over $k$
 such that $\varphi(a v)=a^d\varphi(v)$ for all $a\in k$, $v\in V$ and such that the map $\theta : V \times
\dots \times V \to k$ defined by
 $$\theta(v_1,\dots,v_d)=\frac{1}{d!} \sum_{1\leq i_1< \dots<i_l\leq d}(-1)^{d-l}\varphi(v_{i_1}+ \dots +v_{i_l})$$
($1\leq l\leq d$)  is a $d$-linear form over $k$. The {\it dimension} of $\varphi$ is defined as ${\rm dim}\,\varphi={\rm dim}\, V$.
By fixing a basis $\{e_1,\dots,e_n\}$ of $V$ any form $\varphi:V\to k$ of degree $d$ can be viewed as a homogeneous polynomial
of degree $d$ in $n={\rm dim}\, V$ variables $x_1,\dots,x_n$
 via $\varphi(x_1,\dots,x_n)=\varphi(x_1e_1+\dots+x_ne_n)$.
Conversely, any homogeneous polynomial of degree $d$ over $k$ is a form of degree $d$ over $k$.
Any $d$-linear form $\theta : V \times \dots \times V \to k$ induces a form $\varphi: V\to k$ of degree $d$ via
$\varphi (v)=\theta(v,\dots,v)$.
We can therefore identify $d$-linear forms and forms of degree $d$.
A $d$-linear space $(V,\theta)$ (or the associated $d$-linear form $\theta$) is
called {\it nondegenerate} if $v = 0$ is the only vector such
that $\theta (v, v_{2}, \dots, v_d) = 0$ for all $v_i \in V$.
 A form of degree $d$ is
called {\it nondegenerate} if its associated $d$-linear form is  nondegenerate.
 From now on only nondegenerate forms will be investigated, unless specified otherwise.

\subsection{}
Two $d$-linear spaces $(V_i,\theta_i)$, $i=1,2$, are called {\it isomorphic} (written
$(V_1,\theta_1)\cong (V_2,\theta_2)$ or just $\theta_1\cong\theta_2$) if there exists a bijective linear map
$f:V_1\to V_2$ such that $\theta_2(f(v_1),\dots,f(v_d))=\theta_1(v_1,\dots,v_d)$ for all $v_1,\dots,v_d\in V_1.$

The {\it orthogonal sum} $(V_1,\theta_1)\perp (V_2,\theta_2)$ of two $d$-linear spaces $(V_i,\theta_i)$, $i=1,2$, is defined
to be the $k$-vector space $V_1\oplus V_2$ together with the $d$-linear form
$$(\theta_1 \perp\theta_2)(u_1+v_1,\dots,u_d+v_d)=\theta_1(u_1,\dots,u_d)+\theta_2(v_1,\dots,v_d)$$
($u_i\in V_1$, $v_i\in V_2$). The {\it tensor product} $(V_1,\theta_1)\otimes (V_2,\theta_2)$
 is the $k$-vector space $V_1\otimes V_2$ together with the $d$-linear form
$$(\theta_1 \otimes \theta_2)(u_1\otimes v_1,\dots,u_d\otimes v_d)=\theta_1(u_1,\dots,u_d)\cdot\theta_2
(v_1,\dots,v_d)$$ [H-P].

 A $d$-linear space $(V,\theta)$ is
called {\it decomposable} if $(V,\theta)\cong (V_1,\theta_1)\perp (V_2,\theta_2)$ for two nonzero $d$-linear spaces
$(V_i,\theta_i)$, $i=1,2$. A nonzero $d$-linear space $(V,\theta)$ is called {\it indecomposable} if it is not
decomposable.  $(V, \theta)$ is called {\it absolutely indecomposable} if $(V, \theta) \otimes_k \overline{k}$ is an indecomposable
$d$-linear space over $\overline{k}$, the algebraic closure of $k$.
There exists a Krull-Schmidt Theorem for nondegenerate $d$-linear spaces, if $d \geq 3$: each
 nondegenerate
$d$-linear space  $(V, \theta)$ over $k$ has a decomposition as a direct sum of nondegenerate indecomposable
$d$-linear spaces which is unique up to order and isomorphism [H, 2.3].
 If we can write $\varphi$ in the form $a_1x_1^d+\ldots +a_mx_m^d$, we use the notation $\varphi=\langle  a_1,\ldots
,a_n\rangle  $ and call the form $\varphi$ {\it diagonal}.
 $\varphi$ is called {\it separable} if it becomes isomorphic to a diagonal form over the algebraic closure of $k$.
A form $\varphi$ of degree
$d$ is called \emph{isotropic} if $\varphi(x)=0$ for some nonzero element $x\in V$, otherwise it is called
\emph{anisotropic}.

\subsection{}
Let $l/k$ be a finite field extension and $s:l\to k$ a non-zero $k$-linear map. If $\Gamma:V\times\dots\times V\to
l$ is a nondegenerate $d$-linear form over $l$ then $s\Gamma:V\times\dots\times V\to
k$ is a nondegenerate $d$-linear form over $k$ with $V$ viewed as a $k$-vector space. The $d$-linear
space $(V,s\Gamma)$  is also denoted by $s_*(V,\Gamma)$ or $s_*(\Gamma)$.
If the map $s$ is the trace of a field extension $l/k$, we  simply write $tr_{l/k}(\Gamma)$.
If $\theta$ is a $d$-linear space over $k$ then there exists a canonical isometry
$s_*(\theta_l\otimes_l \Gamma)\cong \theta \otimes_k s_*(\Gamma).$
If  $(V, \Gamma)$ is a nondegenerate $d$-linear space over $l$, then
$s_*(V, \Gamma) $
is a nondegenerate $d$-linear space over $k$ which is indecomposable over $k$ if and only if
$(V, \Gamma)$ is indecomposable over $l$ [Pu].

\subsection{}
 Let $(V,\varphi)$ be a form over $k$ of degree $d$ in $n$ variables over $k$.
 An element $a\in k$ is {\it represented by} $\varphi$ if there is an $x\in V$ such that $\varphi(x)=a$.
 An element $a\in k^{\times}$ such that $\varphi\cong a\varphi$ is called a {\it similarity factor} of the form $\varphi$.
Write $D_k(\varphi)= \{a \in k^\times\mid\varphi(x)=a \text{ for some } x\in V\}$ for the set of non-zero
elements represented by a form $\varphi$ of degree $d$ on $V$ over $k$ and $G_k(\varphi)=\{a\in
k^\times\mid\varphi \cong a\varphi\}$ for the group of similarity factors of $\varphi$ over $k$.
The subscript $k$ is omitted if it is
clear from the context that $\varphi$ is a form over the base field $k$. If $D_k(\varphi)=k^\times$ then
$\varphi$ is called {\it universal}. If $D_l(\varphi)=l^\times$ for each field extension $l$ over $k$,
$\varphi$ is called {\it strongly universal}.  A form $\varphi$  is called \emph{round} if $D(\varphi)\subset G(\varphi)$.

Note that $k^{\times d}\subset G(\varphi)$ and that $G(\varphi)=k^\times$
if $\varphi$ is the zero form. Since $k^{\times d}$ acts as the identity, $G(\varphi)$ is a union of cosets
of $k^\times$ modulo $k^\times/k^{\times d}$. Therefore, we may also think of $G(\varphi)$ as a group of
classes of $d^{th}$ powers; that is, we can work instead with the factor group $G(\varphi)/k^{\times d}$.
By abuse of notation, we sometimes also regard $D(\varphi)$ as a subset of $k^\times/k^{\times d}$ instead of writing
$D(\varphi)/k^{\times d}$.

\subsection{} A nondegenerate form $\varphi(x_1,\ldots,x_n)$ of degree $d$ in $n$ variables
\emph{permits composition} if $\varphi(x)\varphi(y)=\varphi(z)$  where $x$, $y$ are systems of $n$
indeterminates and where each $z_l$ is a bilinear form in $x,y$ with coefficients in $k$. In this case the vector space
$V=k^n$ admits a bilinear map $V\times V\rightarrow V$ which can be viewed as the multiplicative structure of a
(nonassociative) $k$-algebra $A=V$ and so $\varphi(vw)=\varphi(v)\varphi(w)$ for all $v, w\in V$.

 Let $k$ be a field of characteristic $0$ or
$>2d$. Suppose that $\varphi$ is a nondegenerate form of degree $d$ on a finite dimensional unital Jordan
algebra $A$  over $k$ of dimension $n$.  $\varphi$ \emph{permits Jordan composition} if
$\varphi(\{vwv\})=\varphi(v)^2\varphi(w)$ and $\varphi (1)=1$ for all $v, w\in A$ with $\{vwv\}=2(v\cdot w)\cdot
v-w\cdot (v\cdot v)$. If $A$ is a special Jordan algebra, then $\{vwv\}=vwv$ where the
product on the right is given with respect to the associative multiplication in the ambient algebra of $A$.

\subsection{Multiplicative forms}

 Let $\varphi$ be a form over $k$ of degree $d$ and of dimension $n$.

\begin{remark}  (i) If  $\varphi$ is strongly multiplicative then $\varphi$ is multiplicative.\\
(ii) If $D_l(\varphi)$ is a group for every field extension $l$ over $k$ then $\varphi$ is multiplicative.\\
(iii)  If $G_l(\varphi)=l^\times$ for every field extension $l$ over $k$ then $\varphi$ is strongly multiplicative.\\
(iv) If $D_l(\varphi)=G_l(\varphi)$ for every field extension $l$ over $k$ then $\varphi$ is strongly multiplicative.\\
 (v)  If $d=2s+1$ is odd and $\varphi$ is strongly  multiplicative (resp. multiplicative),
 then $a\varphi$ is strongly multiplicative (resp. multiplicative) for any $a\in k^\times$.
\end{remark}

 For  $\varphi$ to be Jordan multiplicative is similar to $\varphi$ permitting Jordan composition,
only that it is required that the $z_l$ are linear depending on the $y_i$'s and quadratically depending on the $x_i$'s.
Thus we may think of the formula
$\varphi(X)^2\cdot\varphi(Y)=\varphi(z_1,\dots,z_n)$
for suitable $z_1,\dots,z_n\in L= k(X,Y)$ as some kind of generalized Jordan composition for $\varphi$.

\begin{remark}  (i) If  $\varphi$ is strongly Jordan multiplicative then $\varphi$ is Jordan multiplicative.\\
 (ii) If  $\varphi$ is strongly multiplicative, then $\varphi$ is strongly Jordan multiplicative.\\
 (iii) If $D_l(\varphi)$ is a group for every field extension $l$ over $k$ then $\varphi$ is Jordan multiplicative.\\
 (iv) If $l^{\times 2}\subset G_l(\varphi)$ or if
 $D_l(\varphi)^{ 2}\subset G_l(\varphi)$ with $D_l(\varphi)^{ 2}=\{a^2\,|\, a\in D_l(\varphi)\}$
   for every field extension $l$ over $k$, then $\varphi$ is strongly Jordan multiplicative.\\
   (v) If $d$ is odd, then $\varphi$ is strongly multiplicative
if and only if $\varphi$ is strongly Jordan multiplicative.\\
 (vi)  If $d=2s$ is even and $\varphi$ is strongly Jordan multiplicative (resp. Jordan multiplicative),
 then $a^s\varphi$ is strongly Jordan multiplicative (resp. Jordan multiplicative) for any $a\in k^\times$.
\end{remark}

\begin{proposition} (i) $\varphi$ is strongly
multiplicative if and only if there exist $z_1,\dots,z_n$ that are linear forms in $Y=(y_1,\dots,y_n)$
with coefficients in $k(X)=k(x_1,\dots,x_n)$ such that
$$\varphi(x_1,\dots,x_n)\cdot \varphi(y_1,\dots,y_n)=\varphi(z_1,\dots,z_n).$$\\
(ii)  $\varphi$ is strongly Jordan
multiplicative if and only if there exist $z_1,\dots,z_n$ that are linear forms in $Y=(y_1,\dots,y_n)$
with coefficients in $k(X)=k(x_1,\dots,x_n)$ such that
$$\varphi(x_1,\dots,x_n)^2\cdot \varphi(y_1,\dots,y_n)=\varphi(z_1,\dots,z_n).$$
\end{proposition}

The proof is analogous to the one of  [L, 2.11] for $d=2$.

 For cubic forms of dimensions $4,6,7,8$ or $\geq 10$ (or quartic forms of dimensions $7,11, 13,14,15$ or $\geq 17$),
there is no composition formula where the $z_i$'s are bilinear forms in $X,Y$ [S, p.140].
If a form permits Jordan composition, then we may choose
 the $z_i$'s as linear forms in $Y$ and quadratic forms in $X$.
  Therefore, the above conclusion to take the $z_i$'s as linear forms in $Y$ might be the strongest one possible.
\begin{remark} Let $\varphi$ be a homogeneous polynomial in $n$ indeterminates over $k$ such that
$\varphi(X)^2 \equiv 1\,{\rm mod}\, K^{\times d}$ for $K=k(x_1,\dots,x_n)$. Then
 $\varphi$  trivially is strongly Jordan multiplicative. 
 For instance, $\varphi(X)=q(X)^{d/2}$ with $q(X)$ a quadratic form in $n$ indeterminates over $k$, is a strongly
Jordan multiplicative form of degree $d$ over $k$.
\end{remark}

\section{Two classes of strongly multiplicative forms}

 Let $A$ be a finite dimensional unital $k$-algebra.
 There exists a nondegenerate form $\varphi$ of degree $d>2$ permitting composition on $A$ if and only if
$A$ is a separable alternative algebra which can be written as a direct sum of simple ideals
$$A=A_1\oplus\dots\oplus A_r$$
 with the center of every $A_i$ a separable field extension $k_i$ of $k$. Any $a\in A$
can be written uniquely as $a=a_1+\ldots +a_r,\, a_i\in A_i$. $\varphi$ is given by
$$\varphi(a)=N_1(a_1)^{s_1}\cdots N_r(a_r)^{s_r},$$
 where $d=d_1s_1+\ldots +d_rs_r$ for positive integers $s_1,\dots,s_r$,
and where $N_i$ is the generic norm of the $k$-algebra $A_i$ (of degree $d_i$).
Viewed as homogeneous polynomials, the norms $N_i$ are irreducible (Schafer [S]).
Forms permitting composition  constitute important examples of nondegenerate
 multiplicative  forms:

\begin{lemma} Let $(A,\varphi)$ be a nondegenerate form of degree $d$ over $k$ permitting composition. Then: \\
(i) $D_l(\varphi)=G_l(\varphi)$ (that is,  $\varphi_l$ is round), for all field extensions $l$ over $k$.
\\
(ii) $\varphi$ strongly multiplicative.
\end{lemma}

\begin{proof} (i) Let $a\in D(\varphi),$ $a=\varphi(y)$. Then
$0\neq y\in A$ is an anisotropic vector and the left multiplication  $L_y:A\rightarrow A,\, x\mapsto yx$ yields
an isomorphism $\varphi\cong a\varphi$, since $a\varphi(x)=\varphi(yx)=\varphi(L_y(x))$ and since $L_y$ is bijective:
suppose that $yx=0$, then $0=\Theta (yx,yx_2,\ldots, yx_d)=\Theta(x,\ldots , x_d)\varphi(y)$ for the $d$-linear
form $\Theta$ associated with $\varphi$, hence $\Theta(x,x_2,\ldots , x_d)=0$ for all $x_2,\ldots , x_d\in A$.
This implies $x=0$, because $\Theta$ is nondegenerate. Therefore $D(\varphi)\subset G(\varphi)$.
The other inclusion is trivial.
For each field extension $l$ over $k$, the form $\varphi\otimes_k l$ will still permit composition.
Therefore $D_l(\varphi)= G_l(\varphi)$.
\\(ii) follows from (i) and Remark 1 (ii).
\end{proof}

 Let $A$ be a finite
dimensional unital commutative algebra over a field $k$ of characteristic $0$ or $>2d$. Then
there exists a nondegenerate form $\varphi$ of degree $d>0$ permitting Jordan composition on $A$ if and only if
$A$ is a separable Jordan algebra which can be written as a direct sum of ideals
$$A=A_1\oplus\ldots\oplus A_r,$$
 where each $A_i$ is a simple Jordan algebra of degree $d_i$ over $k$, whose center is a
separable field extension of $k$, such that $d=d_1s_1+\ldots +d_rs_r$ is satisfied for positive integers
$s_1,\dots,s_r$. Any $a\in A$ can be written uniquely as $a=a_1+\ldots +a_r,\, a_i\in A_i$. $\varphi$ is given by
$$\varphi(a)=N_1(a_1)^{s_1}\cdots N_r(a_r)^{s_r},$$
 where $N_i$ is the generic norm of the $k$-algebra $A_i$ [S].
The generic norm of a Jordan algebra over an arbitrary field is an irreducible homogeneous polynomial [J].
Forms permitting Jordan composition constitute an important class of nondegenerate Jordan
multiplicative forms:

\begin{proposition} Let $(A,\varphi)$ be a nondegenerate form of degree $d$ over $k$
permitting Jordan composition. Then:\\
(i) $D_l(\varphi)^2\subset G_l(\varphi)$ for all field extensions $l$ over $k$, i.e.; $\varphi$ is strongly
Jordan multiplicative.\\
(ii) If $d$ is odd, then $D_l(\varphi)=G_l(\varphi)$ (that is, $\varphi_l$ is round) for all field extensions
 $l$ over $k$, and $\varphi$ is strongly multiplicative.
\end{proposition}

\begin{proof} (i) Let $a\in D(\varphi),$ $a=\varphi(u)$. Then $u\in A$ is  invertible and the
operator $F_u:A\rightarrow A$,
$$F_u(v)=\{uvu\}=2(uv)u-vu^2$$
 is bijective.
  Hence we have
an isomorphism $\varphi\cong a^2 \varphi$, since $a^2\varphi(v)=\varphi(F_u(v))$ for all $v\in A$.
$\varphi$ permits Jordan composition under all field extensions, thus $D_l(\varphi)^2\subset G_l(\varphi)$ for all
extensions $l$ over $k$.
\\
(ii) $1\in D(\varphi)$, hence $G(\varphi)\subset D(\varphi)$.
Let $a\in D(\varphi), a=\varphi(y)$. Then
$a^2\in G(\varphi)$ by (i). Since we also know that trivially $a^{-d}\in G(\varphi)$, this implies that
$a^{-d}a^2 =a^{-(d-2)},a^{-(d-2)}a^2 =a^{-(d-4)},\dots,a^{-(d-2r)} \in G(\varphi)$ for any integer $r>0$,
 since $G(\varphi)$ is multiplicatively closed. Let $d=2r_0+1$ be odd, then $a^{-1}=a^{-(d-2r_0)}\in G(\varphi)$
 and with the same argument also
 $a=a^2a^{-1}\in G(\varphi)$. Hence $D(\varphi)\subset G(\varphi)$.
Since $\varphi$ remains a form permitting Jordan composition under each field extension,
the first assertion is clear. Moreover, $\varphi(X)\in D_K(\varphi)=G_K(\varphi)$, therefore $\varphi$ is strongly
multiplicative.
 \end{proof}

\begin{example}
 Let $k$ be a field of characteristic $0$ or $>6$. Let $\varphi=n_{A/k}$ be
  the norm of a separable simple cubic Jordan algebra $A$.
Then $\varphi$ is absolutely indecomposable [P1].  If $A$ is a reduced Albert algebra then $\varphi$
 $D_l(\varphi)=G_l(\varphi)=l^\times$ for all field extensions
$l$ of $k$. The norm of an Albert division algebra is in general not universal, see [P2]. It is Jordan
multiplicative and a strongly multiplicative cubic form.
 \end{example}

\section{A class of strongly Jordan multiplicative quartic forms}

Let $k$ be a field of characteristic not 2 or 3. To keep the exposition short, we make free use of the terminology
from [A-F1] and refer the reader to [A1, A-F1] for other definitions and results.
An algebra $A$ with involution $\can$ is a unital nonassociatve algebra which is finite-dimensional as a $k$-vector space
together with an anti-automorphism $\can$ of period 2. $A$ is a {\it structurable algebra} if
$$[V_{x,y},V_{z,w}]=V_{V_{x,y}z,w}-V_{z,V_{y,x}w}$$
for all $x,y,z,w\in A$ where $V_{x,y}\in{\rm End}_k(A)$ is defined by
$$V_{x,y}(z)=\{x,y,z\}=(x\bar y)z+(z\bar y)x+(z\bar x)y.$$
 Let $(A,\can)$ be a simple structurable algebra over $k$ such that the space
$S(A,\can)=\{x\in A\,|\, \overline{x}=-x\}$ of skew-symmetric elements has dimension 1. Any such algebra possesses a
nondegenerate quartic form $N_A$ with $N_A(1)=1$ that occurs naturally as the denominator of the inversion operation and induces the
structure of a simple Freudenthal triple system on $(A,\can)$: choose $0\not=s_0\in S(A,\can)$, then there exists an
element $\mu\in
k^\times$ such that $s_0^2=\mu 1$ and define
$$N_A(x)=\frac{1}{12\mu}\chi (s_0x,\{x,s_0x,x\})$$
where $\chi:A\times A\to k$  is a nondegenerate symmetric bilinear form defined by
$$\chi(x,y)=\frac{2}{\mu} \psi(s_0x,y)s_0,$$
 $\psi(x,y)=x\bar y-y\bar x$.
Both $\chi$ and $N_A$ are independent of the choice of $s_0$ [A-F1, p.~190, 191].
An element $x\in A$ is {\it conjugate invertible} (i.e., there exists an element $\hat{u}\in A$ such that $V_{\hat{u},u}=id$)
iff $N_A(x)\not=0$ [A-F1, 2.11]. $(A,\can)$ is called a {\it conjugate division algebra} if every non-zero element of $A$ is conjugate
invertible.

\begin{theorem} Let $(A,\can)$ be a structurable algebra over $k$ such that the space
$S(A,\can)=\{x\in A\,|\, \overline{x}=-x\}$ of skew-symmetric elements has dimension 1.
Let $N_A$ be the nondegenerate quartic form defined above.
Then $D_l(N_A)^2\subset G_l(N_A)$ for all field extensions $l$ over $k$, and $N_A$ is strongly Jordan multiplicative.
\end{theorem}

\begin{proof} Let $a\in D(N_A),$ $a=N_A(u)$. Then $u\in A$ is conjugate invertible and the bijective
operator $P_u:A\rightarrow A$,
$$P_u=U_u(2\tau^{\langle u\rangle}-{\rm id})^{-1}$$
($U_xy=V_{x,y}x$) yields
an isomorphism $N_A\cong a^2N_A$, since $a^2N_A(x)=N_A(P_u(x))$  [A-F1, p.~188, 3.8].
Therefore $D(N_A)^2\subset G(N_A)$.
For each field extension $l$ over $k$, $N_A\otimes_k l=N_{A\otimes l}$.
 Therefore $D_l(N_A)^2\subset G_l(N_A)$.\\
 \end{proof}

\begin{example} ([A-F1]) (i) Let $J$ and $J'$ be $k$-vector spaces possessing
 cubic forms $N$ and $N'$ and paired with
a nondegenerate bilinear form $T:J\times J'\to k $. For $i,j\in J$ and $i',j'\in J'$, define
$j\times i\in J'$ and $j'\times i'\in J$ by
$$T(l,j\times i)=N(j,i,l) \text{ and } T(j'\times i',l')=N'(j',i',l')$$
for $l\in J$, $l'\in J'$ where $N(\,,\,,\,,)$ denotes the trilinear form associated with $N$.
Also, put $j^\sharp=\frac{1}{2} j\times j$ and  $j'^\sharp=\frac{1}{2} j'\times j'$. If the triple satisfies the
{\it adjoint identities}
$$(j^\sharp)^\sharp=N(j)j \text{ and }(j'^\sharp)^\sharp=N'(j')j' $$
and $N$, $N'$ are nonzero, then it is called an {\it admissible triple} ([A1, p.~148], [A-F1, (4.2)]). The $k$-vector space

\smallskip
\[\mathcal{M}(T,N,N')=\{
\left [\begin{array}{ccc}
\alpha &j\\
\noalign{\smallskip}
j'&\beta\\
\end{array}\right ]
\,|\, \alpha,\beta\in k, j\in J,j'\in J'\}\]

\noindent\smallskip
together with the multiplication

\smallskip
\[
\left [\begin{array}{ccc}
\alpha &j\\
\noalign{\smallskip}
j'&\beta\\
\end{array}\right ]
\left [\begin {array}{ccc}
\gamma &i\\
\noalign{\smallskip}
i'&\delta\\
\end{array}\right ]
=
\left [\begin{array}{ccc}
\alpha \gamma + T(j,i')& \alpha i+\delta j+j'\times' i'\\
\noalign{\smallskip}
\gamma j'+\beta i'+j\times i&\beta\delta + T(i,j')\\
\end{array}\right ]
\]

\smallskip\noindent becomes a simple structurable algebra $\mathcal{M}(T,N,N')$
over $k$ such that the space of skew-symmetric elements has dimension 1.
The involution is given by

\smallskip
\[\overline{
\left [\begin{array}{ccc}
\alpha &j\\
\noalign{\smallskip}
j'&\beta\\
\end{array}\right ]
}=
\left [\begin{array}{ccc}
\beta &j\\
\noalign{\smallskip}
j'&\alpha\\
\end{array}\right ].\]

\smallskip\noindent
Furthermore, for

\smallskip
\[x=
\left [\begin{array}{ccc}
\alpha &j\\
\noalign{\smallskip}
j'&\beta\\
\end{array}\right ]
\]
\smallskip\noindent we obtain the strongly Jordan multiplicative quartic form

$$N_A(x)=4\alpha N(j)+4\beta N'(j')-4T({j'}^{\sharp},j^\sharp)+(\alpha\beta -T(j,j'))^2.$$

Let $J$ be a separable Jordan algebra over $k$ of degree 3 with generic norm $N_J$, generic trace $T_J$
and $\zeta\in k^\times$. Then $(\zeta T_J,\zeta N_J,\zeta^2 N_J)$ is a nontrivial admissible triple defined on $(J,J)$.

For every nontrivial triple $(T,N,N')$
satisfying the adjoint identities and defined on a pair of spaces of dimension greater than 2 there are $\zeta$ and
$J$ as above such that $\mathcal{M}(T,N,N')\cong \mathcal{M}(\zeta T_J,\zeta N_J,\zeta^2 N_J)$ (Springer [Sp1, 2]).
If $J$ is a split Albert algebra then we obtain the 56-dimensional irreducible module  for the
split simple Lie algebra of type $E_7$ this way [A-F1, Section 4].
 A structurable algebra $(A,\can)$ of our type is isomorphic to $\mathcal{M}(T,N,N')$ for some triple
 $(T,N,N')$ if and only if $\mu\in k^{\times 2}$ [A-F1, 4.5].
\end{example}

\begin{example} Let $\mu\in k^\times$.
Starting from a separable Jordan algebra $B$ admitting a nondegenerate quartic form $Q$ permitting Jordan composition,
 the $k$-vector space
$${\rm Cay } (B,Q,\mu)=B\oplus B$$
together with the multiplication
$$(b_1,b_2)(b_3,b_4)=(b_1b_3+\mu (b_2b_4^\theta)^\theta, b_1^\theta b_4+(b_2^\theta b_3^\theta)^\theta) $$
becomes a simple structurable algebra whose space of skew-symmetric elements has dimension 1 [A-F1, Section 6].
 Here, $\theta:B\to B$ is a $k$-linear bijection of order 2  defined via
$$b^\theta=-b+\frac{1}{2}t(b)1.$$
The strongly Jordan multiplicative quartic form $N_{{\rm Cay } (B,Q,\mu)}$ is given by
$$N_{{\rm Cay } (B,Q,\mu)}(b_1,b_2)=Q(b_1)+\mu^2Q(b_2)+\frac{\mu}{2} t(U_{b_1}b_2,b_2)-\frac{\mu}{4} t(b_1,b_2)^2,$$
where $t:B\times B\to k$ is the trace form of $Q$.
$(B,id)$ is a subalgebra of $({\rm Cay } (B,Q,\mu),\can)$.

This construction is called the {\it Cayley-Dickson process} (not to be confused with the classical one for quadratic algebras).
If $\mu\in k^{\times 2}$, then ${\rm Cay } (B,Q,\mu)\cong \mathcal{M}(T,N,N')$ for a suitable
triple [A-F1, 6.5].

In particular, let $B$ be a central simple Jordan algebra of degree 4 over $k$ (cf. Section 3), then $B$ has dimension $10,$
$16$ or $28$ [A-F1, 6.12] and ${\rm Cay } (B,Q,\mu)$ has dimension $20$, $32$ or $56$.
For certain choices of $B$, $Q$ and $\mu$, $N_{{\rm Cay } (B,Q,\mu)}$ is anisotropic [A-F1, 7.1].
 We note that $N_{{\rm Cay } (B,Q,\mu)}$ is
not only made up from the scalar $\mu$ and the quartic form $Q$, the algebra structure again plays a role, too.
\end{example}

\begin{remark} (due to B. Allison) Under the hypotheses of [A-F2, Proposition 4.7], $N(P_u(x))=\alpha_uN(P_u(x))$
holds for some  scalar $\alpha_u\in k^\times$   more generally for a structurable algebra whose space of skew-hermitian
 elements may have any dimension. If  $N(\hat{u})=N(u)^{-1}$ holds then $\alpha_u=N(u)^2$ and we would get again
$N(P_u(x))=N(u)^2 N(P_u(x))$. (The norm $N$ may be of different degree than 4 in this case, see [A-F2].)
\end{remark}

\section{How to construct multiplicative forms}

\subsection{} By multiplying different forms we can build classes of absolutely indecomposable (strongly) multiplicative forms of
 higher degree: let $(V_i,\varphi_i)$ be a nondegenerate form of degree $d_i>0$ for every $i$, $1\leq i\leq r$,
 and let $s_1,\ldots , s_r>0$ be integers. Then
 $$(1)\,\,\,\,\,\,\,\,\varphi(u)=\varphi_1(u_1)^{s_1}\cdots \varphi_r(u_r)^{s_r},$$
$u=u_1+\dots +u_r\in V_1\oplus\dots\oplus V_r$ ($u_i\in V_i$ for $1\leq i\leq r$)
 is a form of degree $d=d_1s_1+\ldots +d_rs_r$ on the vector space $V=V_1\oplus\dots\oplus V_r$.
 We can also write $\varphi$
as $$\varphi(x_1,\dots,x_r) =[(\varphi_1\perp 0\dots \perp 0) \cdot
(0\perp \dots \perp\varphi_i \perp\dots \perp 0) \cdot (0\perp \dots \perp 0 \varphi_r)](x_1,\dots,x_r).$$

In the following, suppose as usual that $k$ has characteristic 0 or greater than $d=d_1s_1+\ldots +d_rs_r$.
\begin{remark} (i)  If $a_i\in G(\varphi_i)$, then $a_1^{s_1}\cdots a_r^{s_r}\in G(\varphi)$.\\
(ii) $\varphi$ is anisotropic if and only if $r=1$ (i.e if $\varphi(x)=\varphi_1(x)^{s_1}$)
 and $\varphi_1$ is anisotropic.\\
(iii) If $d\geq 3$, $r\geq 2$ and the $\varphi_i:V_i\to k$ (of degree $d_i\geq 1$) are all nondegenerate, then
$\varphi$ is absolutely indecomposable [Pr, 5.1].\\
(iv) If $\varphi_1$ is any nondegenerate quadratic form over $k$ then the quartic form $\varphi(x)=\varphi_1(x)^2$
is absolutely indecomposable ([K-W, Proposition 10],
in the proof different scalars are used to define the forms).
\end{remark}

\begin{lemma}  Let $(V,\varphi)$ be as in (1).\\
(i) If there is an $i_0$, $1\leq i_0\leq r$ such that $s_{i_0}=1$ and $\varphi_{i_0}$ satisfies $D_l(\varphi_{i_0})=l^\times$
for each field extension $l$ over $k$ (that is, $\varphi_{i_0}$ is strongly universal), then $$D_l(\varphi)=l^\times$$
for each field extension $l$ over $k$ and $\varphi$ is multiplicative and  Jordan multiplicative.
\\(ii) If there is an $i_0$, $1\leq i_0\leq r$ such that $s_{i_0}=1$ and $\varphi_{i_0}$ satisfies $G_l(\varphi_0)=l^\times$
for each field extension $l$ over $k$, then $$G_l(\varphi)=l^\times$$
for each field extension $l$ over $k$ and $\varphi$ is strongly multiplicative.\\
(iii) If every $\varphi_i$ is strongly multiplicative, then so is $\varphi$.\\
(iv) If every $\varphi_i$ is multiplicative, then so is $\varphi$.\\
(v) If every $\varphi_i$ is strongly Jordan multiplicative, then so is $\varphi$.\\
(vi) If every $\varphi_i$ is  Jordan multiplicative, then so is $\varphi$.
\end{lemma}

The proof is trivial and uses Remark 1 (iii) and Remark 2 (iii).

 In particular, if $\varphi_1:V\to k$ is a form of degree $d_1$
 such that $G_l(\varphi_1)=l^\times$ for all field extensions $l$ of $k$, then   $\varphi(u)=\varphi_1(u)^s$
is a strongly multiplicative
 form and $G_l(\varphi)=l^{\times s}$  for all field extensions $l$ of $k$.
If $\varphi_1$ is a form such that $D_l(\varphi_1)=l^\times$ for all field extensions $l$ of $k$,
then $\varphi(u)=\varphi_1(u)^s$
 is a multiplicative form and $D_l(\varphi)=l^{\times s}$  for all field extensions $l$ of $k$. This way one
  can construct
 examples of strongly multiplicative forms of degree $sd_1$ with $k^{\times s}$ as group of similarity factors, and
 examples of multiplicative forms which are not universal.

\begin{example} Let $\varphi$ be as in (1) such that $s_{i_0}=1$ for one $i_0$, $1\leq i_0\leq r$.\\
(i) If $\varphi_{i_0}$ is an isotropic quadratic form over $k$, then $\varphi$
is multiplicative.
\\(ii) If $\varphi_{i_0}$ is either  a hyperbolic quadratic form, $\varphi_{i_0}: k\to k$, $\varphi_{i_0}(u_{i_0})=u_{i_0}$,
or $\varphi_{i_0}$ is the reduced norm  of a split central simple algebra
${\rm Mat}_l(k)$, then $\varphi$
is strongly multiplicative.
\end{example}

\begin{example} Let $k$ be a field of arbitrary characteristic for a change.
Let us  consider the monomials $x_1^{m_1}\cdots x_n^{m_n}\in k[x_1,\dots,x_n]$,
${m_1},\dots, {m_n} >0$. Each such monomial is a form over $k$ in the more
general sense of forms used in [Pr] (see also Roby [R]), which coincides with our definition of forms for fields of
characteristic $0$ or greater than $d$. It is indecomposable unless it is of the kind
\begin{enumerate}
\item $x_1x_2$ and ${\rm char} \, k\not = 2$;
\item $x_1^{p^n}x_2^{p^n}$, where $p\not=2$ is a prime, and ${\rm char}\, k=p$
\end{enumerate}
[Pr, 5.3].
Obviously, (1) is the hyperbolic plane $q\cong \langle 1,-1\rangle$, and as such a strongly multiplicative
quadratic form. Thus over fields of characteristic $0$ or greater than $d$, the hyperbolic plane is the only
monomial which yields a decomposable form of degree $d$.
 Indeed,
$\varphi(x_1,x_2)=x_1^mx_2^m$ must be strongly multiplicative and absolutely indecomposable for each integer $m$ over any field of
characteristic $0$ or greater than $2m$, and $D_l(\varphi)=G_l(\varphi)=k^{\times m}$.

More generally, the form $\varphi:k^r\to k,$
$$\varphi((x_1,\dots,x_r))=x_1^{m_1}\dots x_r^{m_r}$$
of degree $d=m_1+\dots +m_r$  must be strongly multiplicative and absolutely indecomposable for each integer $m$
 over any field of
characteristic $0$ or greater than $d$.
\end{example}

\subsection{}  In the more general setting of algebras over rings, Schm\"{a}hling [Sch] proved
that the composite of two forms of degree $m$ and $n$ (in the sense of [Sch]) with values in an algebra
which both ``permit composition'' is a form of degree $mn$ with values in an algebra ``permitting composition'' [Sch, Satz 2].
 In our setting, this implies nothing new: since we are interested in forms with values in a base field, the forms constructed
this way appear already in Schafer's classification. However, one can generalize this approach:
let $l$ be a finite separable field extension of $k$ with norm $n_{l/k}$. If $(V,\varphi_0)$ is a form of degree $d_0$
and dimension $n$ over $l$, then $\varphi:V\to k$ defined via
$$\varphi(v)=n_{l/k}(\varphi_0(v))$$
is a form of degree $[l:k]\cdot d_0$ over $k$ and of dimension $[l:k]\cdot n$.
We have ${\rm rad} \varphi_0 \subset {\rm rad} \varphi$ where the {\it radical} of an arbitrary form is defined to be
${\rm rad}(\varphi)=\{x\in V\,|\,\theta(x,x_2,\dots,x_d)=0 \text{ for all }x_i\in V\}$ [Sch,  Satz 1, (18)].

\begin{remark} (i) If $\varphi_0$ permits composition, then so does $\varphi(v)=n_{l/k}(\varphi_0(v))$.\\
(ii) If $\varphi_0(u)=\varphi_1(u)^{s_1}\dots\varphi_r(u)^{s_r}$ then $\varphi(u)=
n_{l/k}(\varphi_1(u))^{s_1}\dots n_{l/k}(\varphi_r(u))^{s_r}$.\\
(iii) $\varphi$ is anisotropic if and only if $\varphi_0$ is anisotropic.
\end{remark}

\begin{proposition} Let $l$ be a finite separable field extension of $k$ with norm $n_{l/k}$ and
let $(V,\varphi_0)$ be a form of degree $d_0$ over $l$.
\\(i) If $(V,\varphi_0)$ is strongly multiplicative (respectively, a strongly Jordan multiplicative),
 then $\varphi:V\to k$,
$$\varphi(v)=n_{l/k}(\varphi_0(v))$$
is a strongly multiplicative (respectively, a strongly Jordan multiplicative) form of degree $[l:k]\cdot d_0$ over $k$.\\
(ii) If $(V,\varphi_0)$ is multiplicative (respectively, a Jordan multiplicative), then $\varphi:V\to k$,
$$\varphi(v)=n_{l/k}(\varphi_0(v))$$
is a multiplicative (respectively, a Jordan multiplicative) form of degree $[l:k]\cdot d_0$ over $k$.
\end{proposition}

The proof is a straightforward calculation.

\begin{example} Let $l$ be a finite separable field extension of $k$.
Let $\varphi_0$ be an $r$-fold quadratic Pfister form (or a hyperbolic form of dimension $2m$) over $l$. Then $\varphi(v)=n_{l/k}(\varphi_0(v))$
is a strongly multiplicative form over $k$ of degree $2\cdot [l:k]$ and dimension $2^r \cdot [l:k]$, $r\geq 1$  (or
$2m\cdot [l:k]$, $m\geq 1$ ).

For instance, let $\varphi_0=\langle\langle a_1,\dots,a_r\rangle\rangle$ be an  anisotropic $r$-fold quadratic Pfister form,
 $a_i\in k^\times$.
If $l=k(\sqrt{c})$ is a quadratic field extension, then
$$\begin{array}{l}
n_{l/k}(\varphi_0)(u_1,w_1,\dots,u_{2^r},w_{2^r})=\\
(\langle\langle a_1,\dots,a_r,c\rangle\rangle)^2 (u_1,u_2,\dots,u_{2^r},w_1,w_2,\dots,w_{2^r})-4c\varphi_0^2(u_1w_1,\dots,
u_{2^r}w_{2^r})
\end{array}$$
is an anisotropic strongly multiplicative quartic form of dimension $2^{r+1}$.

If $l=k(\sqrt[3]{c})$ is a cubic Kummer field extension and $k$ contains a primitive third root of unity, then
$$\begin{array}{l}
n_{l/k}(\varphi_0)(u_1,v_1,w_1,\dots,u_{2^r},v_{2^r},w_{2^r})=\\
(\langle\langle a_1,\dots,a_r,2c\rangle\rangle)^3 (u_1,\dots,u_{2^r},v_1w_1,\dots,v_{2^r}w_{2^r})\\
+c(c\langle\langle a_1,\dots,a_r\rangle\rangle\perp 2\langle\langle a_1,\dots,a_r\rangle\rangle)^3
(w_1,\dots,w_{2^r},u_1v_1,\dots,u_{2^r}v_{2^r})\\
+c^2(\langle\langle a_1,\dots,a_r\rangle\rangle\perp 2\langle\langle a_1,\dots,a_r\rangle\rangle)^3
(v_1,\dots,v_{2^r},u_1w_1,\dots,u_{2^r}w_{2^r})\\
-3c[(\langle\langle a_1,\dots,a_r,2c\rangle\rangle(u_1,u_2,\dots,u_{2^r},v_1w_1,\dots,v_{2^r}w_{2^r}))\\
\cdot
((c\langle\langle a_1,\dots,a_r\rangle\rangle\perp 2\langle\langle a_1,\dots,a_r\rangle\rangle)
(w_1,\dots,w_{2^r},u_1v_1,\dots,u_{2^r}v_{2^r}))\\
\cdot
(\langle\langle a_1,\dots,a_r\rangle\rangle\perp 2\langle\langle a_1,\dots,a_r\rangle\rangle)
(v_1,\dots,v_{2^r},u_1w_1,\dots,u_{2^r}w_{2^r}))]
\end{array}$$
is an anisotropic strongly multiplicative  form of degree 6 and dimension $3\cdot 2^{r}$.
\end{example}

 \begin{remark} Let $\varphi_i$, $i=1,2,3$, be forms of degree $d$ over $k$.\\
(i) If $\varphi_1$  satisfies $G_l(\varphi_1)=l^\times$ for all field extensions
 $l$ of $k$ (so, in particular, it is strongly multiplicative), then
 $\varphi=\varphi_1\otimes \varphi_2$ satisfies $G_l(\varphi)=l^\times$ for all field extensions
 $l$ of $k$ and is strongly multiplicative.\\
 (ii) If $\varphi_1$  satisfies $D_l(\varphi_1)=l^\times$ for all field extensions
 $l$ of $k$ (so, in particular, it is  multiplicative), then
 $\varphi=\varphi_1\perp \varphi_2$ satisfies $D_l(\varphi)=l^\times$ for all field extensions
 $l$ of $k$ and is multiplicative.
\end{remark}

In particular, let $\varphi$ be a
form of degree $d$ over $k$ such that $G_l(\varphi)=l^\times$ for all field extensions $l$ of $k$. Then
 $tr_{l/k}(\varphi \otimes l)\cong \varphi \otimes_k  tr_{l/k}(\langle 1\rangle)$ is strongly multiplicative for every separable
field extension $l$ over $k$ and $a_1 \varphi\perp\dots\perp a_n\varphi$ is strongly multiplicative for all $a_1,\dots,a_n\in k^\times$.
This way we are able to construct examples of strongly multiplicative (resp. multiplicative) forms which
  are decomposable by construction and
  irreducible when viewed
 as homogeneous polynomials (e.g., choose $\varphi_1=det$ to be the norm of ${\rm Mat}_d(k)$, or to be the norm of a reduced
 Albert algebra,
  then $\varphi=m\times\varphi_1 =\varphi_1\perp\dots\perp \varphi_1$ is strongly multiplicative).

\section{Restricting the  structure of decomposable strongly multiplicative forms}

The (strongly) multiplicative forms which were constructed as products
$\varphi:V_1\oplus V_2\to k,$
$\varphi(u_1+u_2)=\varphi_1(u_1) \varphi_2(u_2)$ of forms of lower degree are absolutely
indecomposable.

 We now collect some necessary conditions for a decomposable form over a field of characteristic 0 or greater than
 $d$ to be strongly multiplicative.  Let $K=k(X)$, $X=(x_1,\dots,x_n)$ and
 let
$$\varphi_K=\varphi_1\perp\dots\perp\varphi_r$$
 be a decomposable form of dimension $n$ and degree $d\geq 3$ with  the $\varphi_i$'s indecomposable, $r >1$.
For  strongly multiplicative $\varphi$, there exists a permutation $\sigma\in S_r$ such that
$$\varphi(X)\varphi_i\cong  \varphi_{\sigma (i)}$$
over $K$ for all $i\in\{1,\dots,r\}$ [H, 2.3].

\begin{theorem} (i) If $\varphi$ is a strongly multiplicative nondegenerate diagonal form,
then $\varphi\cong\langle 1\rangle$.\\
(ii)  If $\varphi$ is a strongly multiplicative nondegenerate separable form,
then $\varphi\cong\langle 1\rangle$.
\end{theorem}

\begin{proof} (i) Let $\varphi=\langle a_1,\dots,a_n\rangle$ be a nondegenerate diagonal form which is strongly multiplicative.
Then there is a permutation
$\sigma\in S_n$ with $\varphi(X)\langle a_i \rangle\cong\langle a_{\sigma(i)}\rangle$ for all $i\in\{1,\dots,n\}$,
 thus $\frac{a_i}{a_{\sigma(i)}}\varphi(X)\in K^{\times d}$.
  The fact that $\frac{a_i}{a_{\sigma(i)}}\varphi(X)$
  is a diagonal form implies that
   $\frac{a_i}{a_{\sigma(i)}}\varphi(x_1,\dots,x_n)=(c_1x_1+\dots+c_nx_n)^d$,
since $k[x_1,\dots,x_n]$ is a Unique Factorization Domain. Since ${\rm char}\,k>d$
this can only happen if $\varphi(x_1)=(a_1x_1)^d$, that is, if $\varphi(x_1)=x_1^d$.\\
(ii) Let $\varphi$ be a separable form over $k$ which is strongly multiplicative, i.e.
$\varphi(X)\varphi_K\cong \varphi_K$. Let $\overline{k}$ be the algebraic closure of $k$ and put $\overline{K}=\overline{k}(X)$.
Then $\varphi_{\overline{k}}$ becomes a diagonal form and $\varphi(X)\varphi_{\overline{K}}\cong \varphi_{\overline{K}}$.
The generic polynomial $\varphi(X)$ becomes
the generic polynomial $\varphi_{\overline{k}}(X)$ of $\varphi_{\overline{k}}$ which has diagonal form. Using (i), we obtain the assertion.
\end{proof}

\begin{remark} (i) If we drop our general restriction on the characteristic of the base field, we can
find binary diagonal forms which are strongly multiplicative: let ${\rm char}\,k=p$ with $p\not=2$ a prime.
The monomials in Example 5 (2) are all of even degree $d=2p^n$  (hence the characteristic is smaller than $d$ here)
and isomorphic to
the diagonal form $\varphi\cong\langle 1,-1\rangle$ of degree $2p^n$, since $x_1^{p^n}x_2^{p^n}=(z_1+z_2)^{p^n}
(z_1-z_2)^{p^n}$ via the invertible
change of variables $x_1=z_1+z_2$,  $x_2=z_1-z_2$, i.e.,  $$\varphi(x_1,x_2)=q(x_1,x_2)^{p^n}$$
 with $q$ the hyperbolic plane.
These decomposable forms satisfy  $$\varphi(x_1,x_2)^2=q(x_1,x_2)^{2p^n}=q(x_1,x_2)^{d}\equiv 1 \,{\rm mod}\, K^{\times d},$$
 hence
are strongly Jordan multiplicative. They are even strongly multiplicative, since $q$ is.
\\(ii) Let $\varphi(x_1,\dots,x_n)= a_1x_1^d+\dots+a_nx_n^d$, $d=mp$ ($m\geq 1$) be a diagonal form such that
$\varphi(x_1,\dots,x_n)\in K^{\times d}$, where $k$ is a field of odd characteristic $p$. Then it may be that
$\varphi(x_1,\dots,x_n)=(c_1x_1+\dots+c_nx_n)^d$ also for some $n\geq 2$. So in this case - as already observed in (i)
 - there are diagonal forms of dimension greater than 1 which are strongly multiplicative.
\end{remark}

\begin{proposition} (i) If $\varphi$ is a strongly Jordan multiplicative  diagonal cubic
 form, then $\varphi\cong\langle 1\rangle$.\\
 (ii) If $\varphi$ is a strongly Jordan multiplicative  separable cubic
 form, then $\varphi\cong\langle 1\rangle$.
\end{proposition}

\begin{proof} (i) Let $\varphi=\langle a_1,\dots,a_n\rangle$ be a nondegenerate diagonal form which is strongly Jordan multiplicative.
Then there is a permutation
$\sigma\in S_n$ such that $\varphi(X)^2\langle a_i \rangle\cong\langle a_{\sigma(i)}\rangle$ for all $i\in\{1,\dots,n\}$.
 Thus $\frac{a_i}{a_{\sigma(i)}}\varphi(X)^2\in K^{\times 3}$.
 Suppose there is an integer $i_0$ such that $\frac{a_{i_0}}{a_{\sigma(i_0)}}=b^2$ for some $b\in k^\times $.
 Then
 $$\frac{a_{i_0}}{a_{\sigma(i_0)}}\varphi(x_1,\dots,x_n)^2= (ba_1x_1^3+\dots + b a_nx_n^3)^2$$
  is the square of a diagonal form. $(ba_1x_1^3+\dots + b a_nx_n^3)^2\in K^{\times 3}$
 implies that $b\varphi(x_1,\dots,x_n)=(c_1x_1+\dots+c_nx_n)^3$,
since $k[x_1,\dots,x_n]$ is a Unique Factorization Domain (compare the exponents and use that 2 and 3 are coprime).
 It is always assumed that ${\rm char}\,k>3$, therefore this can only happen if $b\varphi(x_1)=(a_1x_1)^3$; i.e.,
 if $b\varphi\cong\langle 1\rangle$, hence $\varphi\cong\langle b^2\rangle$. This, however, means that
$bx^6\langle b^2\rangle\cong \langle b^2\rangle$ or
$bx^6\in K^3$, i.e. $b\in K^3$, thus again $\varphi\cong\langle 1\rangle$.\\(ii)
This is shown  using (i) as in Theorem 2 (ii).
\end{proof}

\begin{lemma}  If $\varphi$ is strongly multiplicative  and  $\varphi_K=\varphi_1\perp\varphi_2$ is its decomposition into
indecomposable components,
then $\varphi(X)\in G_K(\varphi_1)\cap  G_K(\varphi_2)$ or  $\varphi(X)^2\in G_K(\varphi_1)\cap  G_K(\varphi_2)$.
\end{lemma}

Note that if $k$ is a perfect base field, every nondegenerate form of degree $d$ which is indecomposable but not absolutely so,
can be written as the trace of an absolutely indecomposable nondegenerate form over a suitable separable field extension [Pu].

\begin{theorem} Let $l/k$ be a finite separable field extension and $\varphi=tr_{l/k}(\Gamma)$ a
nondegenerate form of degree $d$,
 where $\Gamma$ is an absolutely indecomposable form over $l$. Let $\overline{k}$ be the algebraic closure of $k$
and $\sigma_1,\dots,\sigma_n$ the distinct $k$-algebra isomorphisms from $l$ into $\overline{k}$. If $n<d$ and
$\varphi(X)^s\not\in G(\Gamma_{\overline{K}})$ for all  integers $\,s$, $1\leq s\leq n$, $\overline{K}=\overline{k}(X)$,
then $\varphi$ is  not strongly multiplicative.
\end{theorem}

\begin{proof} Let $l/k$ be a finite separable field extension
and $\sigma_1,\dots,\sigma_n$ the distinct $k$-algebra isomorphisms from $l$ into $\overline{k}$. Then
$\varphi\otimes_k \overline{k}\cong \Gamma^{\sigma_1}_{\overline{k}}\oplus\dots\oplus  \Gamma^{\sigma_n}_{\overline{k}}$
 [H, 4.6]. Suppose $\varphi(X)\varphi_{\overline{K}} \cong\varphi_{\overline{K}}$, then by [H, 2.3],
 there is an integer $s$,  $1\leq s\leq n$, such that $\varphi(X)^s\Gamma_{\overline{K}} \cong\Gamma_{\overline{K}} $,
 which is a contradiction.
\end{proof}

Analogous results hold for strongly Jordan multiplicative forms.

\begin{theorem} Let $K=k(X)$, $X=(x_1,\dots,x_n)$ and let $\varphi$ be a decomposable form of dimension $n$
and degree $d\geq 3$ with indecomposable components $\varphi_i$;
$$\varphi_K=\varphi_1\perp\dots\perp\varphi_r$$
for some $r >1$.\\
(i) If ${\rm dim}\varphi_i\not={\rm dim}\varphi_j$ for $i\not=j$ and there exists an $i_0$ such that
$\varphi(X)\not\in G(\varphi_{i_0})$, then $\varphi$ is not strongly multiplicative.\\
(ii) Suppose there exist $i_0,i_1$ ($i_0\not=i_1$) such that ${\rm dim}\varphi_{i_0}={\rm dim}\varphi_{i_1}$
and ${\rm dim}\varphi_i\not={\rm dim}\varphi_{i_0}$ for all $i\not=i_0$. If there are $h,s\in \{i_0,i_1\}$
such that $\varphi(X)\not\in G(\varphi_h)$ and $\varphi(X)^2\not\in G(\varphi_{s})$, then
$\varphi$ is not strongly multiplicative.  \\
(iii)  Suppose there exist $i_0,i_1$ ($i_0\not=i_1$) with ${\rm dim}\varphi_{i_0}\not={\rm dim}\varphi_{i_1}$
and ${\rm dim}\varphi_i={\rm dim}\varphi_{j}$  as well as
 ${\rm dim}\varphi_{i}\not={\rm dim}\varphi_{i_1},{\rm dim}\varphi_{i_0}$
 for all $i\not=j$, $i,j\not \in \{i_0,i_1\}$.
If there is $h\in \{i_0,i_1\}$ with $\varphi(X)\not\in G(\varphi_h)$, then
$\varphi$ is not strongly multiplicative.
\end{theorem}

The proof is straightforward and uses [H, 2.3]. (An analogous condition can be formulated for
decomposable forms which are not strongly Jordan multiplicative.) Note that for $d$ odd, we have
$\varphi(X)\in G(\varphi_i)$ if and only if $\varphi(X)^2\in G(\varphi_{i})$.

\begin{example} (i) Let $a\in k^\times$ and
 $$\varphi_K\cong\langle a\rangle \perp\varphi_2\perp \dots\perp\varphi_r$$
for some $r >1$, with ${\rm dim}\,\varphi_i >1$ for all $i\geq 2$ and ${\rm dim}\,\varphi_i\not={\rm dim}\,\varphi_j$
for all $i\not=j$. Since
$\varphi(X)\not=X^d$,
 $\varphi$ is not strongly multiplicative (e.g., if
the polynomial representing  one of the $\varphi_i$'s is irreducible).\\
(ii)  Let $a,b\in k^\times$ and
 $$\varphi_K\cong\langle a,b\rangle \perp\varphi_2\perp \dots\perp\varphi_r$$
for some $r >1$, such that ${\rm dim}\,\varphi_i>1$ for all $i >1$.
Since  $\varphi(X)\not=X^d$  and $\varphi(X)^2\not=X^d$, $\varphi$ is not strongly multiplicative.
\end{example}

In analogy to quadratic forms one might ask, for which strongly multiplicative cubic forms $\varphi_0$ and scalars $a\in k^\times$
the cubic form $\varphi_0\perp a\varphi_0$ (respectively, the cubic form $\varphi_0\perp a\varphi_0\perp a^2\varphi_0$)
is strongly multiplicative.

We observe that, although $\varphi_0=\langle 1\rangle$ is strongly multiplicative,
the cubic form $\varphi_0\perp a\varphi_0=\langle 1,a\rangle$ (respectively,  $\varphi_0\perp a\varphi_0\perp a^2\varphi_0
=\langle 1,a,a^2\rangle$) is not strongly multiplicative (Theorem 2), however, the indecomposable norm form $\varphi(u,v,w)=\varphi_0(u)+ a\varphi_0(v)+ a^2\varphi_0(w)
-3auvw$, obtained by the first Tits construction starting with $k$ with scalar $a\in k^\times$, is always strongly
multiplicative.\\
If $\varphi_0$ is strongly multiplicative and $G_K(\varphi_0)=K^\times$, then
 $\varphi_0\perp a\varphi_0$ (respectively,  $\varphi_0\perp a\varphi_0\perp a^2\varphi_0$) is trivially
 strongly multiplicative.
 In particular, if $\varphi_0=n_{A/k}$ is  the norm of a separable simple cubic Jordan algebra $A$ (e.g., a cubic field extension of $k$, a central simple algebra over $k$
or an Albert algebra over $k$), then $\varphi_0$ is strongly multiplicative and
 $\varphi_0\perp a\varphi_0$ (respectively,  $\varphi_0\perp a\varphi_0\perp a^2\varphi_0$) is
 strongly multiplicative for any $a\in k^\times$, provided that $\varphi_{0,K}=n_{A\otimes_k K/K}$
  is surjective.
If $\varphi_{0,K}$ is not surjective and if both $\varphi(X)$ and $\varphi(X)^2$ are not norms of $A\otimes_kK$,
 then the cubic form $\varphi=\varphi_0\perp a\varphi_0$  cannot be strongly multiplicative by Lemma 3.

For cubic forms, the Tits process [P-S1, 2] plays a  role similar to that
of the  Cayley-Dickson doubling.
Hence it is probably not  enough to study cubic forms of the above type.
For instance, let $\varphi_0$ be a cubic form in $n$ variables permitting composition and define
$\varphi(u,v,w)=\varphi_0(u)+ a\varphi_0(v)+ a^2\varphi_0(w)+d\theta_0(u,v,w)$ ($d\in k^\times$) with $\theta_0$ the trilinear form
associated with $\varphi_0$. Here, at least the generic value $\varphi_0(X) $ of the subform
$\varphi_0$ (as well as $a\varphi_0(X) $ or $a^2\varphi_0(X) $) is a similarity factor of $\varphi$
over $K=k(x_1,\dots,x_{n})$.

\section{Some general observations}

We finish this paper generalizing the notion of a multiplicative form of higher degree in two ways.

\subsection{} Assume that the base field $k$ contains a primitive $m$th root of unity $\mu$.

\begin{remark} Let $\varphi$ be a form of degree $d$ over $k$.\\
 (i) If  $\mu \varphi(X) \varphi\cong \varphi$ (resp., $\mu \varphi(X)^2 \varphi\cong \varphi$)  over $K=k(X)$,
 then $\mu \varphi(X) \varphi(Y)\in D_L( \varphi)$ (respectively, $\mu \varphi(X)^2 \varphi(Y) \in D_L( \varphi)$)
 for $L=k(X,Y)$.\\
 (ii) If $\mu \varphi(X) \varphi\cong \varphi$ over $K=k(X)$, then
 $\mu^2 \varphi(X)^2 \varphi\cong \varphi$
over $K=k(X)$.\\
 (iii) If $G_l(\varphi)=l^\times$ for every field extension $l$ over $k$, then $\mu \varphi(X)^2 \varphi\cong \varphi$
over $K=k(X)$.\\
(iv)  $\mu q(X)^2 q\cong q$ over $K=k(X)$ for a quadratic form $q$ if and only if $\mu\in G_k(q)$.
\end{remark}

From now on let $k$ be an infinite field.
For the proof of the next theorem we need the following observation.

\begin{lemma}  Let $f$ and $g$ be polynomials in $k[X]$ such that $f(u_1,\dots,u_n)^m=g(u_1,\dots,u_n)^m$ for all $u_i\in k$, then
there must be an $m$th root of unity $\eta$ such that $f(u_1,\dots,u_n)=\eta g(u_1,\dots,u_n)$ for all $u_i\in k$.
\end{lemma}

\begin{proof}  Let $f$ and $g$ be polynomials in $k[X]$ such that $f(u_1,\dots,u_n)^m=g(u_1,\dots,u_n)^m$ for all
$u_i\in k$, then $f(X)^m-g(X)^m=
\prod_{\mu_i}(f(X)-\mu_i g(X))=0$ (over all $m$th roots of unity $\mu_i$)
in the ring $k[X]$ implying that there must be an $m$th root of unity $\eta$ such that $f(u_1,\dots,u_n)=\eta g(u_1,\dots,u_n)$
 for all $u_i\in k$.
\end{proof}

\begin{theorem} Let $\varphi_1:V\to k$ be a form of degree $>0$ and put $\varphi(v)=\varphi_1(v)^m$, $m>1$.
 $\varphi$ is strongly  multiplicative (resp. strongly Jordan multiplicative) if and only if
$\eta \varphi_1(X) \varphi_1\cong \varphi_1$
(resp. $\eta \varphi_1(X)^2 \varphi_1\cong \varphi_1$) over $K=k(X)$ for some $m$th root of unity $\eta$.
\end{theorem}

\begin{proof} Assume that $\varphi$ is strongly multiplicative. Since $\varphi(X)\varphi(v)
=\varphi(f(v))$
for some isomorphism $f$ and all $v\in V\otimes K$, we get $\varphi_1(X)^m\varphi_1(v)^m=\varphi_1^m(f(v))$ which implies
$\eta\varphi_1(X)\varphi_1(v)=\varphi_1(f(v))$ for some $m$th root of unity $\eta$ by Lemma 4.
 If $\varphi$ is strongly Jordan multiplicative, the proof works analogously.
\end{proof}

\begin{example} Let $\varphi_1:V\to k$ be a form of degree $>0$.
Let $\varphi(v)=\varphi_1(v)^2$. Then $\varphi$ is
\begin{enumerate}
\item strongly multiplicative if and only if $\varphi_1$
is either strongly  multiplicative or  $-\varphi_1(X) \varphi_1\cong \varphi_1$;
\item strongly Jordan multiplicative if and only if $\varphi_1$
is either strongly Jordan multiplicative or  $-\varphi_1(X)^2 \varphi_1\cong \varphi_1$.
\end{enumerate}
If $\varphi_1$ has odd degree $d_1$ (which implies $-1=(-1)^{d_1}$ is a ${d_1}$-th power),
or if $-1\in G_k(\varphi_1)$ then, obviously, $-\varphi_1(X) \varphi_1\cong \varphi_1$
(respectively, $-\varphi_1(X)^2 \varphi_1\cong \varphi_1$) if and only if
$\varphi_1$ is strongly  multiplicative (respectively, strongly Jordan multiplicative).
 \end{example}

\begin{lemma} Let  $\varphi$ be a form of degree $d$.\\
(a) $\mu \varphi(X) \varphi\cong \varphi$ over $K=k(X)$ if and only if
 $\varphi\cong \mu^{d-1} \varphi_0$ with $\varphi_0$ a  strongly multiplicative form.\\
 (b)  $\mu^2 \varphi(X)^2 \varphi\cong \varphi$
over $K=k(X)$ if and only if $\varphi\cong\mu^{d-1}\varphi_0$ with $\varphi_0$ a strongly Jordan multiplicative form.\\
\end{lemma}

\begin{corollary} Let $\varphi_1:V\to k$ be a form of degree $d_1>0$.
Let $\varphi(v)=\varphi_1(v)^m$, $m>1$. Then $\varphi$ is strongly multiplicative if and only if
 $\varphi\cong  \varphi_0^m$ for a strongly multiplicative form  $\varphi_0$.
\end{corollary}

This follows directly from Theorem 5 and Lemma 5. Thus there are no other strongly
multiplicative forms than the obvious ones.

For the remainder of this section we classify  quadratic forms which satisfy $\mu q(X)q\cong q$ over $K=k(X)$.

 Let $ q$ be a nondegenerate quadratic form over $k$.
Obviously, $q$ is (strongly) multiplicative (respectively, strongly multiplicative) if and only if $\mu q$
 satisfies $\mu q(X) q(Y)\in D_L(q)$ for $L=k(X,Y)$ (respectively, $\mu q(X) q\cong q$ over $K=k(X)$).

\begin{theorem} (a) A quadratic form $q$ is isotropic and $\mu q(X)q\cong q$ over $K=k(X)$ if and
only if $q$ is hyperbolic.\\
(b)  For any anisotropic quadratic form $q$ over $k$ the following are equivalent:\\
(i) $q\cong \mu q_1$ with $q_1$ a Pfister form over $k$.\\
(ii) $\mu q(X) q(Y)\in D_L( q)$ for $L=k(X,Y)$.\\
(iii) $\mu q(X) q\cong q$ over $K=k(X)$.
\end{theorem}

\begin{proof} (a) Let $q$ be isotropic such that $\mu q(X)q\cong q$. Then
$ \mu q$ is strongly multiplicative (and isotropic); i.e.,
a hyperbolic form. This yields the assertion.\\
(b) (i) $\Rightarrow$ (iii) is trivial,
 (iii) $\Rightarrow$ (ii) see Remark 8 (i),
 (ii)   $\Rightarrow$ (i) by Lemma 5.
\end{proof}

Therefore an anisotropic  quadratic form such that $\mu q(X) q(Y)\in D_L( q)$ for $L=k(X,Y)$ is either
 an isotropic quadratic form or an anisotropic Pfister form
scaled by $\mu$. A quadratic form such that  $\mu q(X) q\cong q$ over $K=k(X)$ is always a  Pfister form scaled by $\mu$.
In particular, this means:

\begin{corollary}  Let $q$ be a nondegenerate quadratic form and put $\varphi(u)=q(u)^m$, $m>1$.
Then $\varphi$ is strongly multiplicative if and only if $\varphi\cong q_1^m$, $q_1$ a  Pfister form.
\end{corollary}

\begin{lemma} If $\varphi$ is a Jordan multiplicative (respectively, strongly  Jordan multiplicative) form of odd degree $d=2s+1$
over $k$, then $\mu^s\varphi$ satisfies
$\mu \varphi(X)^2 \varphi(Y) \in D_L( \varphi)$ for $L=k(X,Y)$ (respectively, $\mu \varphi(X)^2 \varphi\cong \varphi$
over $K=k(X)$) for any root of unity $\mu$ in $k$.
\end{lemma}

\begin{example} If $\varphi$ is a strongly Jordan multiplicative cubic form, then
$\mu\varphi$ satisfies $\mu \varphi(X)^2 \varphi\cong \varphi$
over $K=k(X)$.
\end{example}

The question is whether there are other forms such that $\mu \varphi(X) \varphi\cong \varphi$ over $K=k(X)$ or
 $\mu \varphi(X) \varphi(Y)\in D_L( \varphi)$ for $L=k(X,Y)$
than the ones constructed above. This has been answered negatively only for quadratic forms so far.
What happens when the form $\varphi_1$ has degree greater than $2$ seems to be much harder to check.
 To find new classes
of such forms could mean we have found new classes of multiplicative form as well.
We conclude with an elementary observation for this more difficult case:
 Let $\varphi_1$ be a form of degree $d_1$.
If $\varphi_1$ is strongly multiplicative (respectively,  strongly Jordan multiplicative) then
 $\varphi=\mu^{d-1}\varphi_1$ satisfies  $\mu \varphi(X)^2 \varphi\cong \varphi$
(respectively, $\mu \varphi(X)^2 \varphi\cong \varphi$) over $K=k(X)$ for some root of unity $\mu$.

\subsection{Strongly multiplicative forms of exponent $s$}

As soon as forms of higher degree than 2 are considered, it makes sense to also define the following notion of multiplicativity:
 Let $\varphi$ be a form over $k$ of degree $d$ and dimension $n$ and let $s$ be a fixed integer, $1\leq s\leq d-1$.
The form $\varphi$ is called {\it strongly multiplicative of exponent $s$} if $\varphi(X)^s\in G_{k(X)}(\varphi)$;
that is, if
$$\varphi(X)^s \varphi\cong \varphi$$
 over $K=k(X)$, where $X=(x_1,\dots,x_n)$.

\begin{remark} (i) Every strongly multiplicative form is  strongly multiplicative of exponent $s$ for every $1\leq s\leq d-1$.
\\(ii) Obviously, $\varphi$ being strongly multiplicative of exponent $2$ is the same as $\varphi$ being strongly
Jordan multiplicative.\\
(iii) If $\varphi_1$ is a form such that $G_l(\varphi_1)=l^\times$ for all field extensions $l$ over $k$,
then $\varphi(u)=\varphi_1(u)^m$ satisfies $l^{\times m}\subset G_l(\varphi_1)$ for all field extensions $l$ over $k$
implying that $\varphi$ is strongly multiplicative of exponent $m$.\\
(iv) Let $e={\rm gcd}\,(s,d)$. Then $\varphi$ is strongly multiplicative of exponent $s$ if and only if
$\varphi$ is strongly multiplicative of exponent $e$. Hence we only have to consider the case $s|d$ when studying
forms which are strongly multiplicative of exponent $s$.
\end{remark}

  Let $\varphi_i$, $i=1,2$, be two forms of degree $d$ over a field $k$
 such that $\varphi_1$ satisfies $l^{\times s}\subset G_l(\varphi_1)$ for all field extensions
 $l$ of $k$ (so  it is, in particular, strongly multiplicative of exponent $s$). We note the following:
 If $\varphi_2$ is a form which represents $1$, then
 $\varphi=\varphi_1\otimes \varphi_2$ satisfies $l^{\times s}\subset
 G_l(\varphi)$ for all field extensions
 $l$ of $k$ and  is strongly multiplicative of exponent $s$. For arbitrary  $\varphi_2$,
 $$\varphi=\varphi_1\otimes \varphi_2\otimes\dots\otimes \varphi_2 \,\, (\text{$s$-times }\varphi_2 )$$
  satisfies $l^{\times s}\subset G_l(\varphi)$ for all field extensions
 $l$ of $k$ and  is strongly multiplicative of exponent $s$.

 \begin{lemma} Let $\varphi$ be any  form of degree $d$ over a field $k$
 such that $\varphi$ is strongly multiplicative of exponent $s$ for some fixed  integer $s$ with $2\leq s\leq d-1$.\\
(i) If $d=r_0s-1$ for a suitable integer $r_0$ then $\varphi$ is strongly multiplicative.\\
(ii) If $d=r_0s+1$ for a suitable integer $r_0$ then $\varphi$ is strongly multiplicative.\\
(iii) Let $r>2$ be an integer such that $rs-d\leq s$. If $\varphi$ is strongly multiplicative of exponent $s$,
then $\varphi$ is strongly multiplicative of exponent $s_0$ with $s_0=r s-d$.
\end{lemma}

\begin{example} (a) Let $d=3$. Each cubic form which is
strongly Jordan multiplicative is strongly multiplicative by (ii). Hence for cubic forms only the notions strongly
 multiplicative and strongly Jordan multiplicative are relevant.
 \\
(b) Let $d=4$. If $\varphi$ is
 strongly multiplicative of exponent $2$, it does not need to be strongly multiplicative.
 \\
(c) Let $d=5$.
 Every form of degree $5$
 which is strongly multiplicative of exponent $s$, $2\leq s\leq 4$, is a strongly multiplicative form by (i), (ii).\\
(d)  Let $d=6$.
 If $\varphi$ is strongly multiplicative of exponent $5$, then it is strongly multiplicative by (ii).
 Furthermore, $\varphi$ is strongly multiplicative of exponent $4$ if and only
 if it is strongly multiplicative of exponent $2$ (one direction is trivial, the other follows with (iii)).
If $\varphi$ is strongly multiplicative of exponent $2$, $3$ or $4$, $\varphi$ does not need to be strongly multiplicative. \\
(e)  Let $d=7$. As before,
 every form of degree $7$ which is strongly multiplicative of exponent $s$, $2\leq s\leq 6$, is a strongly multiplicative form.
\end{example}

The above examples suggest that, at least for all primes up to $7$, it suffices to investigate the strongly multiplicative forms
alone. For $d=4$ or $6$, however, it does make sense to also consider the notion of a form  being
 strongly multiplicative of exponent 2.

\smallskip
{\it Acknowledgements:} The author would like to acknowledge the support of the  the ``Georg-Thieme-Ged\"{a}chtnisstiftung''
(Deutsche Forschungsgemeinschaft) during her stay at the University of Trento,
 and to thank the Department of Mathematics of the University of Trento for its hospitality.

\end{document}